\documentclass[a4paper,12pt]{article}
\usepackage{stmaryrd}
\usepackage{bbm}
\usepackage{amsfonts}
\usepackage{amsmath}
\usepackage{amssymb}
\usepackage{mathrsfs}
\usepackage{accents}
\usepackage{amscd}

\newtheorem{atheorem}{\bf \temp}[section]
\newtheorem{thm}[atheorem]{Theorem}

\newtheorem{lem}[atheorem]{Lemma}

\newtheorem{de}[atheorem]{Definition}

\numberwithin{equation}{section}

\title{\textbf{$L^{4,\infty}$-solution to the Navier-Stokes Equations in four-dimensional space}}
\author{  Xixia Ma \footnote{Corresponding author. E-mail addresses: kfmaxixia@163.com(Xixia Ma)} \ \ \ \ \ \ \
 \\}
\date{}

\begin{document}

\maketitle

\textbf{Abstract.} In this paper, we have showed $L^{4,\infty}$-solutions of the cauchy problem for the four-dimensional Navier-Stokes equations are smooth through backward uniqueness and analytic function's properties .


\begin{center}
\item\section{Introduction}
\end{center}

In this paper ,we consider the Cauchy problem for incompressible Navier-Stokes equations in four spatial dimension with unite one viscosity and zero external force :
\begin{equation}
\begin{array}{l}
 \partial_{t}u-\triangle u+u\nabla u+\nabla p=0\\
 div\quad u=0
\end{array}
\end{equation}
for $x\in \mathbb{R}^{4},t>0$, together with the initial condition \begin{equation} u(0,x)=a(x) \end{equation} in $ x\in \mathbb{R}^{4}  .$

Many authors have studied the regularity of the solutions of the Navier-Stokes equations ,especially spatial dimension $ d=3$ . Ladyzhenskaya-Prodi-Serrin have given a criterion about the regularity of the Navier-Stokes equations, later Struwe,M has showed that the critical points are also right if $ d\geq3$ ,except $ L^{\infty}$ in time direction .As the suitable weak solution was introduced by Caffarelli-Kohn-Nirenberg , they have proved that there is an open subset that $ u $ is H$\ddot{o}$lder continuous and the $1-D $ Hausdorff measure of the singular set equals zero , if $ d=3$and for any suitable weak solution $(u,p)$.Lin.F.H gave a more direct proof of Caffarelli,Kohn,Nirenberg's proof .Hongjie Deng,Dapeng Du gave the proof of four-dimension solution of the Navier-Stokes equations  ,overcoming the difficulties of the solution's compactness .Their method was to apply the backward heat kernel to the Navier-Stokes equations . By the uniqueness of the backward heat kernel ,L.Esacauriaza ,G.Seregin , V.$\check{S}$ver$\acute{a}$k gave a proof of the critical points $ L^{\infty}(L^{3})$ if $ d=3 $. In this paper , we will show that at the critical points $ L^{\infty;}(L^{4})$ if $ d=4$ , the conclusion is also right through the backward uniqueness and lemma 2.2.Furthmore this conclusion shows that the critical condition implies it that the solution of the Navier-Stokes equations is small.However, we want to know whether the conclusion is right or not if $ d\geq5$ .

\begin{thm}
suppose that $ a\in J^{\circ}$
hold. Let $ v $ and $ p $ be a weak Leray-Hopf solutions to the Cauchy problem $(1.1)$and$(1.2)$. Assume that .for some $T>0$the velocity field $v$ satisfies the so-called Ladyzhenskaya-Prodi-Serrin condition ,i.e.
$v\in L^{s,l}(Q_{T})$ with
$$\frac{4}{s}+\frac{2}{l}=1, s\in(4,\infty].$$ 
Then $ v $ is locally essential bounded in  $ Q_{T}$
where $ Q_{T}=\mathbb{R}^{4}\times(0,T)$.\end{thm}

(The proof is found in Struwe,M."On partial regularity results for the Navier-Stokes equations.Comm.Pure Appl.Math.41(1988),no.4,437-458") A Leray-Hopf weak solution of the Cauchy problem(1.1)and (1.2) in  $ Q_{T}$ is a vector field $ v $ satisfying:

a. $$  \quad\quad v\in L^{2,\infty}(Q)\cap L^{2}(-1,0;H^{1}(B)) .     $$

b. the function
$$ t\rightarrow \int_{\mathbb{R}^{4}}v(x,t)w(x,t)dx $$ is continuous in$[0,T]$ for
any $ w\in L^{2}. $

c.$ v $ satisfies the Navier-Stokes equations in the distribution sense and the local energy inequality is valid .

 We recall that the norm in the mixed Lebesgue space $L^{s,l}(Q_{T})$ is given as follows:\begin{equation}
       \|f\|_{s,l;Q_{T}}=\begin{cases}
    (\int\|f(\cdot,t)\|^{l}_{s}dt)^{\frac{1}{l}}, \quad \quad l\in[1,+\infty) \\ esssup_{t\in(0,T)}\|f(\cdot,t)\|_{s},l=\infty
     \end{cases}
     \end{equation}

   In this paper we address the problem of regularity for the weak Leray-Hopf solutions $v$ satisfying the additional condition $v\in L^{4,\infty}(Q_{T})       $.
      we prove that the Theorem 1.1 is valid .More precisely ,we have
  \begin{thm}Assume that $v$ is a weak Leray-Hopf solution to the Cauchy problem $(1.1)$and $(1.2)$in $Q_{T}$ and satisfies
   the additional condition $v\in L^{4,\infty}(Q_{T})       $.
   Then $v\in L^{6,6}(Q_{T})   $
   and hence it is smooth  in $Q_{T}$ .\end{thm}

 In fact ,we prove the following local result .

    \begin{thm}Consider two functions $v$ and $p$ defined in the space-time cylinder $Q=B\times(0,1)$, where $B(r)\subseteq R^{4}$ stands for the ball of radius $r$ with the center at the origin and $B=B(1)$. Assume that $v$ and $p$ satisfy the Navier-Stokes equations in $Q$ in the sense of distributions and have the following differentiability properties:
   $$ v\in L^{2,\infty}(Q)\cap L^{2}(-1,0;H^{1}(B)) ,  p\in L^{\frac{3}{2}}(Q)    $$
  In addition 
                    $$\|v\|_{4,\infty}( Q)<\infty           $$
   Then $v$ is H$\ddot{o}$lder continuous in the closure of the set
            $ Q(\frac{1}{2})=B(\frac{1}{2})\times(-(\frac{1}{2})^{2},0)$.\end{thm}

\begin{center}
\item\section{Suitable Weak Solutions and Backward Uniqueness}
\end{center}

 In this section ,we are going to discuss  smoothness of the so-called suitable weak
solutions to the Navier-Stokes equations .

 \begin{de}let $\Omega$ be a open set in $\mathbb{R}^{4}$. We say that a pair $u$ and $p$ is a suitable weak solution to the Navier-Stokes equations on the set $\Omega\times(-T_{1},T)$ if it satisfies the conditions :
   \begin{equation}
   u\in L^{2,\infty}(\Omega\times(-T_{1},T))\cap L^{2}(-1,0;H^{1}(\Omega)) ,  p\in L^{\frac{3}{2}}(\Omega\times(-T_{1},T))
    \end{equation}
   $u$ and $p$ satisfy the Navier-Stokes equations in the  distribution sense ;

   $u$ and $p$ satisfy the local energy inequality
   \begin{equation}
    \int_{\Omega}\varphi|u(x,t)|^{2}+2\int_{\Omega\times(-T_{1},t)}\varphi|\nabla u|^{2}dxdt'
    \leq$$ $$\int_{\Omega\times(-T_{1},t)}(|u^{2}(\triangle\varphi+\partial_{t}\varphi)+u\cdot\nabla\varphi(|u|^{2}+\nabla2p))dxdt' \end{equation}
 for a.a. $t\in(-T_{1},T)$ and for all nonnegative functions $\varphi\in C^{\infty}_{0}(\mathbb{R}^{4}) $ vanishing in the neighborhood of  the parabolic boundary$ \partial'Q\equiv\Omega\times{
 {t=-T_{1}}}\bigcup\partial\Omega\times[-T_{1},T]$.\end{de}

   \begin{lem} There exist absolute positive constants $ \varepsilon_{0} $and $c_{0k},k=1,2,\cdot\cdot\cdot$, with the following property .Assume that the pair $ u $ and $ p $ is suitable weak solution to the Navier-Stokes equations in $ Q $ and satisfies the condition
   \begin{equation}
            \int_{Q}(|u|^{3}+|p|^{\frac{3}{2}}dz)<\varepsilon_{0}                         \end{equation}
  Then for any natural number $k$, $\nabla^{k-1}u$ is holder continuous in $\overline{Q}(\frac{1}{2})$ and the following bound is valid :
            \begin{equation}
            max_{z\in{Q}(\frac{1}{2})}u<c_{0k}
            \end{equation} \end{lem}

 \textbf{proof:} $k=0,$ the conclusion is ok though interpolation inequality and iteration ,also see[13] .By Theorem 1.1,
  $k\geq 1$ is also right .

In the following content , we are going to give known facts from the theory of unique continuation for differential inequalities .We also work with the backward heat operator $\partial_{t}+\triangle u$.In the space-time cylinder $Q(r,T)\equiv B(r)\times(0,T)\subseteq \mathbb{R}^{4}\times \mathbb{R}^{1}$,we consider a vector-valued function $u=(u_{i})=(u_{1},\cdot\cdot\cdot,u_{n})$, satisfying three conditions :
          \begin{equation} u\in W^{2,1}_{2} (Q(r,T),\mathbb{R}^{n})    ;    \end{equation}
 \begin{equation} |\partial_{t} u -\triangle u|\leq c_{1}(|u|+|\nabla u|)
\quad a.e. in \quad Q(r,T)    \end{equation}
for some positive constant $c_{1}$ ;
 \begin{equation}
   |u(x,t)|\leq C_{k}(|x|+\sqrt{t})^{k}    \end{equation}
for all $k=0,1,\cdot\cdot\cdot,$ for all $(x,t)\in Q(r,T)$, and for some positive constants $C_{k}$. Here
$$ W^{2,1}_{2} (Q(r,T)\equiv{|u|+|\nabla u|+|\nabla^{2}u|+|\partial_{t}u|\in L^{2}(Q(r,T))}.$$

\begin{thm} Assume that a function $u$ satisfies conditions $(2.5)-(2.7)$. Then $u(x,0)=0$ for all $x\in B(r)$\end{thm}(the proof is founded in the part IV of [1])

\begin{thm}  Let $\mathbb{R}^{n}_{+}={{x|(x_{i}\in R^{n}),x_{n}>0},i=1,2,...,}$ and $Q_{+}=\mathbb{R}^{n}_{+}\times(0,1)$,$u:Q_{+}\rightarrow \mathbb{R}^{n}$, and satisfies
 \begin{equation}|\partial_{t}+\triangle u|\leq c_{1}(|u|+|\nabla u|) \quad in \quad Q_{+}   .  \end{equation}For some $c_{1}>0$ and \begin{equation}u(x,0)=0 \quad in \quad R^{n}_{+}       .  \end{equation} We also assume $u\in W^{2,1}_{2}(Q'_{+})$ where $Q'_{+}\subseteq Q_{+} is\quad bounded $. Then if \begin{equation}|u(x,t)\lfloor\leq e^{M|x|^{2}},                                        \end{equation}we have  $ u\equiv0 $  in $ Q_{+}$. \end{thm}

   \begin{thm}Let $ u $ be a solution of (1),(2) such that $ u(\cdot,t)$ is analytic in a bounded open set $ Q=\Omega\times(0,T)$,If there exist a nonempty open set $ \Omega_{1}$ in $ \Omega $ and a $ t_{1}\in(0,T) $ such that $ u(x,t_{1})=0, x\in \Omega_{1} $,then $ u\equiv 0 $ in $ Q $.\end{thm}

 \textbf{proof}:\quad Since $ u(x,t)$ is analytic in   $ x $  and  $ t $ in $Q$ .By assumption $ u(x,t_{1})=0,$ for $ x\in \Omega_{1}$,hence $ u(x,t_{1})=0$ in $\Omega $ .So $\omega(x,t_{1})=0 $  in  $\Omega $. Since $\omega $ satisfies $$\partial_{t}\omega-\triangle \omega=\ast d (u\nabla u)=div (u\wedge \omega).$$ We have $\partial_{t}\omega(x,t_{1})=0 $ and so $\ast du_{t}(x,t_{1})=0 $ .Since $ u_{t}\in H^{1}_{0}$ and $ div  u_{t} =0 $ we deduce  $ u_{t}(x,t_{1})=0 $.Applying the same argument ,we have $ \frac{\partial}{\partial t}^{k}u(x,t_{1})=0 $ for $ k=0,1,2,...$, then the theorem is proved .

\begin{center}.
\item\section{Proof of the main results}
\end{center}

  We start with the proof of theorem $1.3$and then state that Theorem 1.2 is valid .

  Our approach is based on the theory of analytic function's property and backward uniqueness for the heat operator .We use the backward uniqueness to prove Theorem 1.3.Before the proof of the theorem 1.3 ,we give a lemma .

  \begin{lem} Assume $(v,p)$ is a Leray-Hopf solution of (1.1) and (1.2), if $ v\in L^{4,\infty}(Q_{T})$, then we have $\partial_{t}v,\nabla p,\nabla^{2} v\in L^{\frac{4}{3}}(Q_{T})$ .\end{lem}

   \textbf{proof}: First by H$\ddot{o}$lder inequality and Sobolev embedding,
   Choose $\phi\in C^{\infty}_{0}(\mathbb{R}^{4})$ and $ div \phi =0$, we have
   \begin {equation}
   \begin{array}{l}
  (\partial_{t}v,\phi)=-(v\nabla v,\phi)-(\nabla v,\nabla\phi)\\
                     \leq\|v\|_{L^{2}}\|\nabla v\|_{L^{2}}\|\phi\|_{L^{\infty}}+\|\nabla v\|_{L^{2}}\|\nabla\phi\|_{L^{2}}\\
                   \leq(\|v\|_{L^{2}}+\|v\|_{L^{2}}\|\nabla v\|_{L^{2}})\|\phi\|_{H^{2}} \end{array} \end{equation}.

                   Hence $\partial_{t}v\in L^{\frac{4}{3}}(Q_{T})$. In the following , we show $\nabla p\in L^{\frac{4}{3}}(Q_{T})$ ,then it is easy to check $\nabla^{2} v\in L^{\frac{4}{3}}(Q_{T}$.

                   In fact , let $ f=\partial_{t}v-\triangle v $ ,then first it is obtained that
                   $f\in L^{2}(0,T;H^{-2}_{0})  $ as mentioned above.

                   And then we know
\begin {equation}
   \begin{array}{l}
    div f =0 \\
   \ast df= \ast d(v\nabla v) \end{array}\end {equation} in any open set $\Omega\subseteq\mathbb{R}^{4}$ for a.e $ t\in(0,T)$ .

   By the elliptic regularity theory , \begin{equation}
   \|f\|^{L\frac{4}{3}}_{L\frac{4}{3}}  \leq\|v\nabla v\|^{L\frac{4}{3}}_{L\frac{4}{3}}+ \|f\|^{L\frac{4}{3}}_{H^{-2}_{0}}.\end{equation}
So we get
$\nabla p\in L^{\frac{4}{3}}(Q_{T})$ .

  \textbf{the proof of Theorem 1.3}  Here we will prove in two steps .
  
  Step 1: we show that $(v,p)$ is a suitable weak solution in $Q$. In fact , we need to prove $\|v\|_{L^{3}{(Q(\frac{3}{4}))}}<\infty$, hence it is enough to show  \begin{equation}sup_{-(\frac{3}{4})^{2}\leq t\leq0}\|v(.,t)\|_{L^{3}(B(\frac{3}{4}))}\leq\|v(x,t)\|_{L^{3}(Q)} . \end{equation}

    Claim:  \begin{equation}\begin{array}{l}t\rightarrow \int_{B(\frac{3}{4})}v(x,t)w(x,t)dx\quad
   is\quad  continuous\quad in\quad  [-(\frac{3}{4})^{2},0]\\ \quad for\quad
any\quad w\in L^{\frac{3}{2 }}(B(\frac{3}{4})). \end{array} \end{equation}

   First ,we know $\|v\|_{L^{4}(Q)}<\infty$ by $ v\in L^{4,\infty}(Q)$,hence
   \begin{equation}\|v\nabla v\|_{L^{\frac{4}{3}}(Q(\frac{3}{4}))}\leq\|v\|_{L^{4}(Q(\frac{3}{4}))}\|\nabla v\|_{L^{2}(Q(\frac{3}{4}))} . \end{equation}
  It is easy to show  \begin{equation}\|\partial_{t} v \|_{L^{\frac{4}{3}}(Q (\frac{3}{4}))}< \infty  \end{equation} by (3.1).And this implies (3.4). We also consider  \begin{equation}\|p\|_{L^{\frac{3}{2}}{(Q(\frac{3}{4}))}}<\infty .\end{equation}Compose  \begin{equation}p=p_{1}+p_{2} \end{equation} where \begin{equation}\triangle p_{1}=div(v\nabla v) \quad in \quad B(\frac{3}{4}),  \end{equation} and $p_{1}=0$ on $\partial B(\frac{3}{4})$,(3.8)is proved by C-Z estimate and the properties of the harmonic function .

  Step 2: Assume that the statement of Theorem 1.3 is false .Let $z_{0}\in\overline{Q(\frac{1}{2})}$ be a singular point,then by Lemma 2.2 ,there exists a sequence of positive numbers $r_{k}$ such that $r_{k}\rightarrow 0$ as $ k\rightarrow \infty$ , and   \begin{equation}A(r_{k})=\int_{B(x_{0},r_{k})\times{t_{0}-r^{2}_{k}\leq t\leq t_{0}}}|v(x,t)^{3}dxdt>\varepsilon_{\ast}
 \quad  for all\quad k\in\mathbb{N} .\end{equation}  Here $\varepsilon_{\ast}$ is an absolute positive constant .

   we extend functions $(v,p)$ to the whole space $\mathbb{R}^{4+1}$ by zero .Extended functions will be denoted by $(\widetilde{v},\widetilde{p})$, respectively .Now ,we let
   $$v^{r^{k}}(x,t)=r^{k}\widetilde{v}(x_{0}+r^{k}x,t_{0}+r^{2k}t),$$$$p^{r^{k}}(x,t)=r^{2k}\widetilde{p}(x_{0}+r^{k}x,t_{0}+r^{2k}t.$$
   $$\phi(x,t)=r^{k}\phi^{r^{k}}(x_{0}+r^{k}x,t_{0}+r^{2k}t$$ where $\phi\in C^{\infty}_{0}(R^{4+1)}$ .

   we choose $r^{k}$ so small to ensure $$spt \phi \subset{(x,t)|t_{0}+r^{2k}t\in(-(\frac{4}{3}^{2},(\frac{4}{3}^{2}),x_{0}+r^{k}x\in B(\frac{3}{4})},$$  $$spt\phi^{r^{k}}\subset B(\frac{3}{4})\times(-(\frac{4}{3}^{2},(\frac{4}{3}^{2}).$$
   Then we have  $$2\int_{B\times(-1,0)}\phi^{r^{k}}|\nabla v|^{2}dxdt
    \leq\int_{B\times(-1,0)}(|v^{2}(\triangle\phi^{r^{k}}+\partial_{t}\phi^{r^{k}})+v\cdot\nabla\phi^{r^{k}}(|v|^{2}+\nabla2p))dxdt,$$
   Like estimating $p$ as above mentioned ,we further obtain  \begin{equation}\int_{{Q_{1}}}(|p^{r^{k}}|^{\frac{3}{2}}+|\nabla v^{r^{k}}(x,t)|^{2})dz \leq c_{3}(Q_{1})<\infty        \end{equation} where $ Q_{1}\subseteq \mathbb{ R}^{4+1}$ with $ c_{3}(Q_{1})$ independent of $r^{k}$.Then we apply Step 1and Lemma3.1 ,we find  \begin{equation}\nabla^{2} v^{r^{k}},\nabla p^{r^{k}}\in L^{\frac{4}{3}}(Q)  . \end{equation}  Together with $v\in L^{4,\infty}(Q)$, this implies  \begin{equation}v^{r^{k}}\rightarrow u\quad in \quad L^{3}(Q_{1})          \end{equation} for $ Q_{1}\subseteq \mathbb{ R}^{4+1}$ ;

    Moreover  by Arzela-Ascoli theorem, it is clear that \begin{equation}v^{r^{k}}\rightarrow u \quad in \quad C([a,b];L^{\frac{4}{3}}(\Omega)) \end{equation}for any $-\infty<a<b<+\infty$ and for any $\Omega\subseteq\mathbb{ R}^{4}   $ .

     So furthermore,\begin{equation}v^{r^{k}}\rightarrow u\quad in\quad C([a,b];L^{2}(\Omega))\end{equation}which is easily obtained from the interpolation inequality .

     Now we conclude that :
   \begin{equation}\int_{Q}(|u|^{4}+|\nabla u|^{2}+|\partial_{t}u|^{\frac{4}{3}}+|\nabla^{2}u|^{\frac{4}{3}}+|\nabla u|^{\frac{4}{3}})dz \leq c_{3}(Q)\end{equation}  for any $ Q\subseteq \mathbb{R}^{4+1}  $ ;and \begin{equation}u\in C([a,b];L^{2}(\Omega))\end{equation} for $-\infty<a<b<+\infty$ and for any $\Omega\subseteq\mathbb{ R} ^{4}  $ , and $(u,p)$ satisfies the Navier-Stokes equations a.e in the distribution sense .

    Hence it is easy to show that $(u,p)$ is a suitable weak solution to the Navier-Stokes equations in $\Omega\times[a,b]$ . According to (3.11),$$ sup_{-r^{2k}\leq t\leq0}\frac{1}{r_{k}}\int_{B(0,r_{k})}|v(x,t)|^{2}dx=sup_{-1\leq t\leq0}\int_{B(0,1)}|v^{r^{k}}(x,t)|^{2}dx>\varepsilon_{\star} $$ by the interpolation inequality .For all $k\in\mathbb{N}$ , and by (3.16), we obtain \begin{equation} sup_{-1\leq t\leq0}\int_{B(0,1)}|u(x,t)|^{2}dx>\varepsilon_{\star}     .\end{equation}

     In the following , we are going to show that there exist some positive numbers$ R_{2}$ and $T_{2}$ such that for any $k=0,1,...,$ the function $\nabla^{k}u $ is h$\ddot{o}$lder continuous and bounded on the set $(\mathbb{R}^{4}\setminus\overline{B}(R_{2}))\times (-2T_{2},0]$.

      Let us fix an arbitrary number $ T_{2}>2 $ and note that
   $$\int_{\mathbb{R}^{4}\times(-T_{2},0)}|u(x,t)|^{3}+|q(x,t)|^{\frac{3}{2}}dxdt<\infty .$$

  Therefore ,
   $$ \int_{(\mathbb{R}^{4}\setminus\overline{B}(R_{2}))\times(-T_{2},0)}|u(x,t)|^{3}+|q(x,t)|^{\frac{3}{2}}dxdt\rightarrow 0 $$ as $ R\rightarrow\infty $.

    This means that there exits a number $ R_{2}(\varepsilon_{0},T)>4 $ such that  \begin{equation}\int_{(\mathbb{R}^{4}\setminus\overline{B}(R))\times(-T_{2},0)}|u(x,t)^{3}+|q(x,t)^{\frac{3}{2}}dxdt<\varepsilon_{0}   .\end{equation}

    Now assume that$  z_{1}=(x_{1},t_{1})\in(\mathbb{R}^{4}\setminus\overline{B}(R_{2}))\times (-2T_{2},0]$ . Then $$ Q(z_{1},1)=B(x_{1})\times (t_{1}-1,t_{1})\subseteq\mathbb{R}^{4}\setminus\overline{B}(\frac{R_{2}}{4}))\times (-4T_{2},0] .$$ So by(3.20),  \begin{equation}\int_{B(x_{1})\times (t_{1}-1,t_{1})}|u|^{3}+|q|^{\frac{3}{2}}<\varepsilon_{0}    \end{equation}
   for any $z_{1}\in(\mathbb{R}^{4}\setminus\overline{B}(R_{2}))\times (-2T_{2},0]$, where $ T_{2}$ and $R_{2}>4$. Then it follows from (3.21) and Lemma 2.2 that for any $k=0,1,...,$ \begin{equation} max_{z\in Q(z_{1},\frac{1}{2})}|\nabla^{k}u(z)|\leq c_{0k}<\infty   .\end{equation}  
and $\nabla^{k}u(z) $ is H
$\ddot{o}$lder continuous on $(\mathbb{R}^{4}\setminus\overline{B}(R_{2}))\times (-2T_{2},0] $.
    Now  let us introduce the vorticity $ \omega $ of $ u $, i.e. $\omega=\ast d u$.The function $\omega$ meets the equation $$\partial_{t}\omega-\triangle \omega=\ast d (u\nabla u)=div(u\wedge \omega)$$ in $(\mathbb{R}^{4}\setminus\overline{B}(R_{2}))\times (-T_{2},0]$.

     Recalling (3.22) the function $\omega$ satisfies the following relations : \begin{equation}|\partial_{t}\omega-\triangle \omega|\leq M(|\omega|+|\nabla\omega|)        \end{equation} for some constant $ M>0$ and   \begin{equation}|\omega|\leq c_{00}+c_{01}<\infty        .  \end{equation}

     Let us show that  \begin{equation}\omega(x,0)=0,    x\in\mathbb{R}^{4}\setminus\overline{B}(R_{2}).       \end{equation} We take into account the fact that $ u\in C([a,b];L^{2})$ and find $$(\int_{B(x_{\ast,1)}}|u(x,0)|^{2}dx)^{\frac{1}{2}}\leq $$
    $$(\int_{B(x_{\ast,1)}}|u(x,0)-v^{r^{k}}(x,0)|^{2}dx)^{\frac{1}{2}}+(\int_{B(x_{\ast,1)}}|v^{r^{k}}(x,0)|^{2}dx)^{\frac{1}{2}}$$
    $$\leq\|u-v^{r^{k}}\|_{C([a,b];L^{2})} + |B|^{\frac{1}{6}}(\int_{B(x_{\ast,1)}}|v^{r^{k}}(x,0)|^{3}dx)^{\frac{1}{3}}.$$

   By (3.4)and (3.16), it is concluded that $$(\int_{B(x_{\ast,1)}}|u(x,0)|^{2}dx)=0$$ for all $ x\in\mathbb{R}^{4}$,so(3.25)is proved .Hence by Theorem 2.4 of Section 2, we show that  \begin{equation}\omega(z)=0,   z\in(\mathbb{R}^{4}\setminus\overline{B}(R_{2}))\times (-T_{2},0]      . \end{equation}

    Claim: \begin{equation}\omega(\cdot,t)=0 \quad in\quad \mathbb{R}^{4}     \quad   for \quad a.e.\quad t\in(-T_{2},0) .\end{equation}

     We know that $(u,q)$ meet the equations :  \begin{equation}\partial_{t}u+u\nabla u+\nabla q=0,divu=0, \triangle u=0, \nabla\wedge u=0       \end{equation} in the set $(\mathbb{R}^{4}\setminus\overline{B}(R_{2}))\times (-T_{2},0]$. From (3.28), we deduce the following bound  \begin{equation} max_{Q_{0}}(|\nabla^{k}u|+|\nabla^{k}\partial_{t}u|+|\nabla^{k}q|\leq c^{1}_{0k}<\infty    \end{equation} for all\quad $k=0,1,2,...$ ,here $Q_{0}=(\mathbb{R}^{4}\setminus\overline{B}(R_{2}))\times (-T_{2},0]$.

      To prove (3.27),according to (3.26) , we fix a smooth cut-off function $\varphi\in C^{\infty}_{0}(\mathbb{R}^{4})$ subjected to the conditions :$\varphi(x)=1$ if $ x\in B(2R_{2})$ and $\varphi(x)=0 $ if $ x\in\mathbb{R}^{4}\setminus B(3R_{2})$ . Let $ w=\varphi u,r=\varphi q $ ,so $(w,r)$ satisfies
       \begin{equation}
       \begin{array}{l}
       \partial_{t}w-\triangle w+w\nabla w+\nabla r=g\\
       div w=u\nabla\varphi \end{array}  \end{equation} in $ Q_{\ast}=B(2R_{2})\times (\frac{-T_{2}}{2},0)$and \begin{equation} w=0 on\quad\partial B(2R_{2})\times (\frac{-T_{2}}{2},0)                  \end{equation} where $ g=(\varphi^{2}-\varphi)u\nabla u+uu\cdot\nabla\varphi^{2}+q\nabla\varphi-2\nabla u\nabla\varphi-u\triangle\varphi$

       It is clear that $w$ is not incompressible .So we introduce the functions $(\widetilde{w},\widetilde{r})$ satisfies :$$-\triangle\widetilde{w}+\nabla\widetilde{r}=0, div\widetilde{w}= u\nabla\varphi $$ in $ Q_{\ast}$ with $ \widetilde{w}=0$ on $\partial B(2R_{2})\times (\frac{-T_{2}}{2},0)$.

        Setting $ U= w-\widetilde{ w} $ and $ P=r-\widetilde{r}$ satisfies  \begin{equation}
        \begin{array}{l}
        \partial_{t}U-\triangle U+U\nabla U+\nabla P=G-div(U\otimes\widetilde{w}+ \widetilde{w}\otimes U)\\
        div U=0\end{array} \end{equation} in $Q_{\ast}$,  and  \begin{equation} U=0  \quad on \quad  \partial B(2R_{2})\times (\frac{-T_{2}}{2},0)    \end{equation}where $ G=-div\widetilde{w}\otimes\widetilde{w}+g-\partial_{t}\widetilde{w} $.By (3.22)and the elliptic regularity theory, we can choose $ t_{0}\in(\frac{-T_{2}}{2},0)$ so that  \begin{equation}\|\nabla U(\cdot,t_{0}\|_{2, B(2R_{2})}<\infty .\end{equation}

         Then by the short time unique solvability results for the Navier-Stokes equations ,we find a number $\delta_{0}>0 $ such that $$\partial_{t}U,\nabla P,\nabla^{2}U \in L^{2}(B(3R_{2})\times(t_{0},t_{0}+\delta_{0})),$$ then it is easy to check that $$sup_{t_{0}-\varepsilon<t<t_{0}+\delta_{0}-\varepsilon}sup_{x\in B(2R_{2})}|\nabla^{k}U|\leq c^{5}_{0k}<\infty $$ for $k=0,1,...$, and for $0<\varepsilon<\frac{\delta_{0}}{4}$, so it is valid that $$sup_{t_{0}+\varepsilon<t<t_{0}+\delta_{0}-\varepsilon}sup_{x\in B(2R_{2})}|\nabla^{k}u|\leq c^{6}_{0k}<\infty .$$

          Hence $ u(\cdot,t)$ is analytic in the $ B(3R_{2})$ for $(t_{0}+\varepsilon,t_{0}+\delta_{0}-\varepsilon)$,and as mentioned above,$ \omega=0 $ for  $ (B(3R_{2})\setminus B(2R_{2}))\times(t_{0}+\varepsilon,t_{0}+\delta_{0}-\varepsilon) $.By the Theorem 2.5 ,we obtain $$\omega=0$$ in $ B(3R_{2})\times(t_{0}+\varepsilon,t_{0}+\delta_{0}-\varepsilon)$. Then theorem 1.3 is proved .

           \textbf{the proof of Theorem 1.2}:we can choose $w=|v|^{2}$ in the proof of Theorem 1.3. By Theorem 1.3 , it is easy to find $$ w \in L^{2,\infty}(Q_{T})\bigcap L^{2}(\delta,T;W^{1}_{2}(\mathbb{R}^{4})),$$ and then $$\|w\|_{L^{3}}\leq\|w\|_{L^{2}}^{\frac{1}{3}}\|\nabla w\|_{L^{2}}^{\frac{2}{3}},$$ so we deduce $$ w\in L^{3}(Q_{\delta,T})\Leftrightarrow v\in L^{6}(Q_{\delta,T})$$ for some $\delta>0 .\Box $
           \begin{center}
\item\section{Appendix}
\end{center}

    Now we give the existence of the solution of the Navier-Stokes equations in Theorem 1.2 .The method is similar to the  construction of Lei-Lin solution.

   Let us recall that the incompressible Navier-Stokes equations in $\mathbb{R}^{+}\times\mathbb{R}^{4}$ are :
  \begin{equation}
\begin{array}{l}
 \partial_{t}u-\triangle u+u\nabla u+\nabla p=0\\
 div\quad u=0
\end{array}
\end{equation}
one assumes that the initial data $ u(x,0)=u_{0}(0)$ are divergence-free and possess certain regularity .

The known a priori Leray-Hopf energy estimates satisfied by classical solutions of (4.1) is as follows :\begin{equation}sup_{t>0}\|u(t,\cdot)\|_{L^{2}}\leq\|u_{0}\|_{L^{2}},\int^{\infty}_{0}\|\nabla u(t,\cdot)\|^{2}_{L^{2}}dt\leq\frac{1}{2}\|u_{0}\|^{2}_{L^{2}} .\end{equation}

Denote $ D=\sqrt{-\triangle}$. Let $ u  $ be a divergence-free vector field .We make the following decomposition:
\begin{equation} u=u_{+}+u_{-} \end{equation}
where $$ u_{+}=\frac{1}{2}(u+D^{-1}\nabla\times u)$$
and  $$ u_{-}=\frac{1}{2}(u-D^{-1}\nabla\times u)$$

We have the following facts :

\textbf{a:} Let $ u\in H^{1}(\mathbb{R}^{n})$ be a divergence-free vector field and be decomposed into $ u_{+}$ and $ u_{-}$ as (4.3) . Then the following identities hold : $$\nabla\times u_{+}=D u_{+} ,\nabla\times u_{-}=-D u_{-}.$$

\textbf{b:} Let $ m,k\geq 0 $ be any integers and $ u\in C^{m}([0,T),H^{k}(\mathbb{R}^{n}))$. Suppose that for each $ t\in[0,T) $, $ u(t,\cdot)$ is divergence-free .Decompose $  u(t,\cdot)$ into $ u_{+}$ and $ u_{-}$ as (4.3).Then for all integers $ m_{1},m_{2}$ and $ k_{1},k_{2}$ with $ m_{1}+m_{2}\leq m $ and $  k_{1}+k_{2}\leq k $, we have $$\int D^{m_{1}}\partial^{k_{1}}_{t}u_{+}\cdot D^{m_{2}}\partial^{k_{2}}_{t}u_{-}\equiv 0 .$$ 

Now we give a example satisfying the conditions of Theorem 1.2 .

For $ j\geq 1 $ ,let $ N_{j}\geq 1 $ be an integer which is increasing with respect to $ j $ ,$ 0 <\delta_{j}<\frac{1}{2}$ is a constant and $ \gamma_{j}$ is a measurable ,bounded function such that \begin{equation}
\begin{array}{l}
 supp \gamma_{j}\subset C_{N_{j}}={\xi\in\mathbb{R}^{4}:||\xi|-N_{j}|\leq\delta_{j}N_{j}},\\
 N_{j}\leq \frac{N_{k}}{8} for j < k
\end{array}
\end{equation}
We define $ v_{0}$ and $ v_{0j}$ as follows :
  \begin{equation} v_{0}=\Sigma^{\infty}_{j=1} v_{0j}, \widehat{ v_{0j}}(\xi)=(N(\xi)+i|\xi|^{-1}\xi\times N(\xi)\gamma_{j}(\xi) , \end {equation}
where $  N(\xi)$ is a unit vector which is perpendicular to $\xi $ .Clearly ,one has
\begin{equation} \nabla\cdot v_{0j}=0 , v_{0j+}= v_{0j}, v_{0j-}=0. \end {equation} 

Define \begin{equation} v(t,x)=\Sigma_{j\geq 1}v_{j}(t,x) \end{equation}
where $ v_{j}=\Sigma_{j\geq 1}e^{t\triangle} v_{0j}$ . Then it is easy to check that $ v $ satisfies :\begin {equation}
\partial_{t}v +v\cdot\nabla v +\nabla p_{1} =\triangle v + g, \nabla \cdot v = 0 \end {equation}
where $ p_{1}=-\frac{1}{2}|v|^{2},g=-\Sigma_{j,k\geq 1}v_{j}\times \omega_{j}$ and $\omega_{j}=\nabla\times v_{j}$.

If we write $ u=v+b $ where $ v $ is defined above . We can show $sup_{t\geq 0}\|u(t,\cdot)\|_{\dot{H}^{1}}<\infty $ ,here the critical energy $ E_{c}(u)=\frac{1}{2}\|Du(t,\cdot)\|^{2}_{L^{2}}+\int^{t}_{0}\|D\nabla u(t,\cdot)\|^{2}_{L^{2}}ds$ .By the Sobolev embedding theorem ,this estimate implies that $$ u\in L^{q}_{t}(L^{p}),\frac{4}{p}+\frac{2}{q}=1,4\leq p\leq \infty .$$
See the details in [15].

I appreciate my advisor professor Binglong Chen and am thank professor Zhen Lei . At the same time , I am glad to discuss with Dr Xuming Gu .

\end{document}